\documentclass[11pt]{amsart}

\usepackage{amsmath,amssymb,amsthm}
\usepackage{graphicx}
\usepackage{subfigure}
\usepackage{enumerate}
\usepackage{fullpage}
\usepackage{url}
\usepackage{float}
\usepackage{xspace}
\usepackage[hidelinks]{hyperref}
\usepackage{color}
\usepackage{thmtools, thm-restate}
\usepackage[numbers]{natbib}
\usepackage{cleveref}
\usepackage{booktabs}
\usepackage{tikz}
\usetikzlibrary{arrows.meta,calc}

\makeatletter
\newcommand{\addresshere}{%
  \enddoc@text\let\enddoc@text\relax
}
\makeatother

\newtheorem{theorem}{Theorem}

\newtheorem{lem}{Lemma}[section]

\newtheorem{remark}{Remark}

\newcommand{\norm}[1]{\lVert #1\rVert}

\title{Upper bounds for packings of cylinders\\ around the unit sphere}

\author{Alexey Glazyrin}

\address{Alexey Glazyrin, School of Mathematical \& Statistical Sciences, The University of Texas Rio Grande Valley, Brownsville, TX 78520, USA}
\email{alexey.glazyrin@utrgv.edu}

\date{}

\begin{document}

\maketitle
\begin{abstract}
In this note we study the problem of packing infinite circular cylinders tangent to the unit sphere. For a given number of cylinders $N$, $N\geq 3$, the maximal common cylinder radius $R_N$ is the maximum $R$ for which there exist $N$ circular cylinders of radius $R$ with non-overlapping interiors, all tangent to the unit sphere. In 1990, Kuperberg \cite{K90} conjectured that there does not exist a packing of more than 6 cylinders of radius one tangent to the unit sphere, that is, $R_7<1$. This conjecture was recently confirmed by an exhaustive computer-assisted case analysis \cite{MR26}. Here we prove new bounds for $R_N$ for small $N$ by considering an auxiliary configuration of spherical quadrilaterals. In particular, we provide a short geometric proof of Kuperberg's conjecture and show that $R_7<0.9978$.
\end{abstract}

\section{Introduction}

By a packing of cylinders in $\mathbb{R}^3$ we mean a finite family of infinite cylinders with non-overlapping interiors. We are interested in packings with the largest possible common radius of circular cylinders tangent to the same unit sphere. The existence of the maximal $R_N$ for a given $N$ is a consequence of compactness of a configuration space. For a fixed radius $r$, each cylinder of radius $r$ tangent to the unit sphere is specified by a point of the projectivized tangent bundle $T=\mathbb P(T\mathbb S^2)$, which is compact. Given $r$, the set of admissible configurations of $N$ non-overlapping cylinders is compact in $T^N$. Let $R$ be the supremum of the radii admitting a configuration. $R<\infty$ since $N\geq 3$. Then taking an increasing sequence $r_k\to R$, we get a nested family of compact sets of admissible configurations in $T^N$. It must have a common point and nonoverlap is preserved in the limit so $R=R_N$ has a nonempty admissible set as well. We also note that for any $r<R_N$, there is a non-overlapping configuration of $N$ cylinders of radius $r$ because we can substitute each cylinder of radius $R_N$ by its subset, a cylinder of radius $r$.

In 1990, Kuperberg \cite{K90} posed the following problem: how many non-overlapping infinite unit cylinders can touch a unit ball? Six parallel cylinders arranged around an equatorial hexagon show that six are possible. The question was whether seven are possible, or, in the notation above, whether $R_7\geq1$. Kuperberg conjectured that the answer is indeed 6 (see also \cite[Section~2.4, Conjecture 16]{BMP05}). He also asked for the values of $R_N$ (see \cite[Section~2.4, Problem~17]{BMP05}).

The first upper bounds were geometric. Heppes and Szab\'o \cite{HS91} proved that at most eight unit cylinders can touch the ball. They also obtained $R_6\leq1.275$, $R_7\leq1.075$, and $R_8\leq0.96$ (see \cite[Theorem~1.2]{Y19}). Bra\ss\ and Wenk \cite{BW00} reduced the unit-radius upper bound to seven in 2000. For each cylinder they consider its intersection with a larger concentric sphere. The resulting patches have disjoint interiors, so their total area cannot exceed the area of that sphere. Varying its radius gives upper bounds for every $N$ (see Firsching \cite[Table~3.2]{F15} for concrete evaluations for $3\leq N\leq11$). For example, the resulting seven-cylinder
bound is approximately $1.069485$. Starostin \cite{S06} extended related area arguments to tubes with curved axes, a broader class of bodies.

The elementary lower bound comes from parallel axes through the vertices
of a regular $N$-gon:
\begin{equation}\label{eq:parallel-lower}
 R_N\geq p_N,\qquad
 p_N=\frac{\sin(\pi/N)}{1-\sin(\pi/N)}.
\end{equation}
It gives $p_3=3+2\sqrt3$ and $p_4=1+\sqrt2$. The three-cylinder value is exact because any packing of cylinders of radius $r$ induces a packing of balls of radius $r$ tangent to the central ball. For three balls, the optimal packing is when their centers form an equilateral triangle in a plane through the origin. The parallel construction need not be optimal for larger numbers of cylinders. Firsching's numerical constructions, first reported in 2014 and described in his dissertation \cite[Section~3.1.3 and Table~3.1]{F15}, give $R_5>1.45289$, $R_6>1.049659$, and $R_7>0.846934$.

Ogievetsky and Shlosman \cite{OS21a} improved the six-cylinder lower bound to
\[
 R_6\geq\frac{3+\sqrt{33}}8=1.0930703308\ldots
\]
by an explicit deformation with $\mathbb D_3$ symmetry. They subsequently proved that this configuration is a strict local maximum modulo simultaneous rotations \cite{OS21b}; its global optimality is still conjectural. Their deformation also enlarges the parallel ring for every even $N\geq6$ \cite[Section~6]{OS21a}. They also developed other constructions of packings of cylinders \cite{OS21c}, albeit not improving the lower bounds used in our table. The survey \cite{OS19} places these constructions and their rigidity properties in a common framework.

Yard\i mc\i\ \cite{Y19} strengthened the area argument by accounting for uncovered portions of the auxiliary sphere. His upper bounds for $R_9$, $R_{10}$, $R_{11}$ are, after upward rounding, $0.754942$, $0.662272$, and $0.591479$. His further gap analysis, developed with Bezdek, excludes seven cylinders of radius $1.0496594$ \cite[Theorem~5.1]{Y19}.

Recently, Mati\'c and Radoi\v ci\'c \cite{MR26} settled the original question by a finite exact-arithmetic proof excluding seven unit cylinders and thus showed that $R_7<1$. Mati\'c, Radoi\v ci\'c, and Spomer \cite{MRS26} proved $R_4=1+\sqrt2$. They also found exact certificates for $R_5<39/25$, $R_6<1203/845$, and $R_8<14558721/18995711$; the six-cylinder bound is weaker than Heppes--Szab\'o's. More generally, they showed that the insertion of one cylinder into the parallel $m$-ring, with rearrangement allowed, is impossible for $3\leq m\leq7$ and possible for every $m\geq9$.
In particular,
\begin{equation}\label{eq:insertion-lower}
 R_N\geq p_{N-1}\qquad(N\geq10).
\end{equation}
Insertion into the eight-cylinder ring remains open. The exact values of $R_N$ remain unknown for all $N\geq5$.

For general $N$, the parallel construction gives $p_N\sim\pi/N$. The concentric-sphere estimate yields $R_N=O(N^{-2/3})$. These general bounds leave a substantial gap. It is not even known whether $NR_N$ is bounded.

In this paper we concentrate on explicit upper bounds for small numbers of cylinders. The following theorem records the cases in which the resulting certified bound improves the earlier results.

\begin{theorem}[Improved upper bounds]\label{thm:bounds}
$R_6<1.252$, $R_7<0.9978$, $R_9<0.726$, $R_{10}<0.642$, $R_{11}<0.577$, $R_{12}<0.526$, $R_{13}<0.483$, $R_{14}<0.448$, $R_{15}<0.418$, $R_{16}<0.392$, $R_{17}<0.370$, $R_{18}<0.350$, and $R_{19}<0.3327$.
\end{theorem}


\begin{remark}
With our approach using auxiliary quadrilaterals, we could not improve upper bounds obtained by \cite[Theorem 3.1]{Y19} for $20\leq N\leq 24$ but it is quite possible that a more accurate choice of auxiliary regions improves bounds for large $N$ as well. 
\end{remark}

Table~\ref{tab:comparison} assembles the current comparison for every $3\leq N\leq19$, with sources in both earlier-bound columns. Bold entries in the last column are the improvements in
Theorem~\ref{thm:bounds}. The $N=4,5,8$ quadrilateral bounds are included for completeness. For $N=3$ the exact answer is elementary, and no separate quadrilateral certificate is needed.

\begin{table}[h]
\centering\small
\setlength{\tabcolsep}{6pt}
\renewcommand{\arraystretch}{1.15}
\begin{tabular}{@{}rllr@{}}
\toprule
$N$&\shortstack[l]{Current lower bound\\and source}
&\shortstack[l]{Previous upper bound\\and source}
&\shortstack[r]{Quadrilateral\\upper bound}\\
\midrule
3&$3+2\sqrt3$ (parallel)&$3+2\sqrt3$&---\\
4&$1+\sqrt2$ (parallel)&$1+\sqrt2$ \cite{MRS26}&$2.686$\\
5&$1.45289^{a}$ \cite{F15}&$39/25$ \cite{MRS26}&$1.705$\\
6&$(3+\sqrt{33})/8$ \cite{OS21a}&$1.275$ \cite{HS91,Y19}&$\mathbf{1.252}$\\
7&$0.846934^{a}$ \cite{F15}&$1$ \cite{MR26}&$\mathbf{0.9978}$\\
8&$0.67149^{b}$ \cite{OS21a}&$0.766422$ \cite{MRS26}&$0.838$\\
9&$0.519803$ (parallel)&$0.754942$ \cite{Y19}&$\mathbf{0.726}$\\
10&$0.519803$ \cite{MRS26}&$0.662272$ \cite{Y19}&$\mathbf{0.642}$\\
11&$0.447213$ \cite{MRS26}&$0.591479$ \cite{Y19}&$\mathbf{0.577}$\\
12&$0.392239$ \cite{MRS26}&$0.535590^{c}$ \cite{Y19}&$\mathbf{0.526}$\\
13&$0.349198$ \cite{MRS26}&$0.490302^{c}$ \cite{Y19}&$\mathbf{0.483}$\\
14&$0.314605$ \cite{MRS26}&$0.452831^{c}$ \cite{Y19}&$\mathbf{0.448}$\\
15&$0.286208$ \cite{MRS26}&$0.421285^{c}$ \cite{Y19}&$\mathbf{0.418}$\\
16&$0.262485$ \cite{MRS26}&$0.394343^{c}$ \cite{Y19}&$\mathbf{0.392}$\\
17&$0.242375$ \cite{MRS26}&$0.371046^{c}$ \cite{Y19}&$\mathbf{0.370}$\\
18&$0.225114$ \cite{MRS26}&$0.350689^{c}$ \cite{Y19}&$\mathbf{0.350}$\\
19&$0.210138$ \cite{MRS26}&$0.332737^{c}$ \cite{Y19}&$\mathbf{0.3327}$\\
\bottomrule
\end{tabular}
\caption{Bounds for $R_N$.
Parallel refers to \eqref{eq:parallel-lower}; the lower bounds from
\cite{MRS26} are \eqref{eq:insertion-lower}, rounded down.
Every entry in the last column is certified and gives a strict upper
bound.}
\label{tab:comparison}
\par\medskip
\begin{minipage}{\textwidth}\footnotesize

\noindent
$^{a}$Firsching's reported numerical constructions for $N=5,7$.

\noindent
$^{b}$Evaluation of the family from \cite{OS21a}, not numerical values printed in that source.

\noindent
$^{c}$These are exact, conservatively rounded evaluations made here of the earlier inequality \cite[Theorem~3.1]{Y19}, not numerical values printed in that source. The earlier upper decimals for $N=8,9,10,11$ are rounded upward. For $N=8$ the underlying exact bound is $14558721/18995711$.
\end{minipage}
\end{table}

\newpage
\section{Main construction}

Fix $r>0$ and consider infinite circular cylinders of radius $r$. For a cylinder $C$ tangent to the unit sphere $\mathbb{S}^2$, we construct a compact spherical region $Q(C)\subset \mathbb{S}^2$ that is a closed topological disk. We require $Q$ to be equivariant under the action of the isometries of the sphere, that is, for any isometry $g$ of $\mathbb{R}^3$ preserving the unit sphere, the region $Q(g(C))$ is the same as $g(Q(C))$. In particular, regions should be symmetric under the symmetries of the cylinder preserving the unit sphere. Note that $|Q(C)|$, the standard spherical area of the region $Q(C)$, is the same for all $C$ so we use the notation $|Q|$ for this area. As a regularity condition, we also require that the boundary of $Q(C)$ is piecewise smooth and the boundaries of $Q(C_1)$ and $Q(C_2)$ have finitely many intersection components for two cylinders $C_1$ and $C_2$.

For a non-negative parameter $\delta$, we define the family $\mathcal{Q}_{r,\delta}$ of equivariant assignments $C\to Q(C)$ such that $Q\in \mathcal{Q}_{r,\delta}$ if it satisfies the following three properties.

\begin{enumerate}

\item\label{it:1} No triple intersections: $Q(C_1)\cap Q(C_2)\cap Q(C_3)=\emptyset$ for any three non-overlapping cylinders $C_1$, $C_2$, $C_3$.

\item\label{it:2} Connected differences: $Q(C_1)\setminus Q(C_2)$ is connected if $C_1$ and $C_2$ are non-overlapping cylinders.

\item\label{it:3} Small pairwise intersections: $|Q(C_1)\cap Q(C_2)|\leq \delta$ for any  non-overlapping cylinders $C_1$ and $C_2$.

\end{enumerate}

Suitable regions can provide a certificate of nonexistence for cylinder packings.

\begin{theorem}\label{thm:main}
Let $Q\in\mathcal{Q}_{r,\delta}$. If for a positive integer $N$, $N\geq 3$,
\[
N |Q|-(3N-6)\delta>4\pi,
\]
then $R_N<r$.
\end{theorem}

\begin{proof}
Assume $R_N\geq r$ so there exist $N$ non-overlapping cylinders of radius $r$ tangent to the unit sphere. For each cylinder we construct a corresponding region and analyze the configuration of auxiliary regions $Q_i$, $1\leq i\leq N$, on the sphere. By property (\ref{it:1}), the nerve formed by the interiors of the regions in the configuration is just a graph.

We now explain that property (\ref{it:2}) implies that this graph is planar. Indeed, put $Q_i'=Q_i\setminus \bigcup\limits_{j\neq i} Q_j$. If some $Q_j$ overlaps with $Q_i$, then it meets the boundary of $Q_i$ in one or more intervals and carves off several regions from $Q_i$ while keeping $Q_i\setminus Q_j$ connected by property (\ref{it:2}). We choose one of these connected regions $A_{ij}$ with nonempty interior arbitrarily but coordinate our choice so that $A_{ji}$, the intersection region we choose for $Q_j$, equals $A_{ij}$. Note that portions on the boundary of $Q_i'$ formed  by boundaries of different $Q_j$ do not intersect because of property (\ref{it:1}).

Removing portions of all $Q_j$ overlapping with $Q_i$ from $Q_i$ leaves a connected region $Q_i'$ with nonempty interior. Inside each $Q_i'$ choose a point $P_i$. $A_{ij}$ has a part of its boundary from the boundary of $Q_j$ in the interior of $Q_i$. We choose a point $P_{ij}$ there. $A_{ji}$ also has a part of its boundary from the boundary of $Q_i$ in the interior of $Q_j$. We choose a point $P_{ji}$ there. Then we connect $P_i$ to all chosen $P_{ij}$ by arcs whose interiors lie in $Q_i'$ and are pairwise disjoint. Finally, we connect every $P_{ij}$ and $P_{ji}$ by an arc whose interior lies in the interior of $A_{ij}$. This way we realize the graph of intersections of the interiors of regions as a geometric graph in the sphere with no edge crossings. Therefore, the graph is planar.

Due to Euler's formula, the number of edges in the graph of intersections is not larger than $3N-6$. We can now estimate the total area of $\bigcup\limits_{i=1}^N Q_i$:
\[
\left|\bigcup\limits_{i=1}^N Q_i\right| = \sum\limits_{i=1}^N |Q_i| - \sum\limits_{1\leq i<j\leq N} |Q_i\cap Q_j|\geq N |Q|  -(3N-6)\delta>4\pi.
\]
The total area exceeds the area of the sphere, which is a contradiction, so this configuration is impossible.
\end{proof}

\begin{remark}
The constructions of Heppes-Szab\'o \cite{HS91} and Bra\ss-Wenk \cite{BW00} are special cases of \Cref{thm:main} for $\delta=0$ after radial projections on $\mathbb{S}^2$. The approach of Yard\i mc\i{} \cite{Y19} is a modification of the $\delta=0$ case: \cite[Theorem 3.1]{Y19} shows that for certain regions with pairwise disjoint interiors, their complement necessarily contains two disjoint spherical caps.

We could not find a suitable region with $\delta=0$ for the proof of Kuperberg's conjecture. We conjecture that it does not exist, that is, there is no $Q\in\mathcal{Q}_{1,0}$ with $|Q|>\frac{4\pi}{7}$.
\end{remark}

\section{Spherical quadrilaterals}\label{sect:quad}
In this section we describe the regions that produced new bounds. We fix a cylinder $C$ of radius $r$. Let its axis be vertical and pass through $(0,-1-r,0)$. Then for $u,v>0$, we define
\begin{equation*}\label{eq:Q}
 Q_{u,v}(C)=\{(x,y,z)\in\mathbb{S}^2:y\leq0,\ |x|\leq u(-y),\ |z|\leq v(-y)\}.
\end{equation*}
This is a geodesically convex quadrilateral in the open hemisphere $y<0$. Its center is the cylinder's tangency point $c=(0,-1,0)$ and its vertices are
\begin{equation}\label{eq:vertices}
 V_{\sigma\tau}=\frac{(\sigma u,-1,\tau v)}{\sqrt{1+u^2+v^2}},
 \qquad \sigma,\tau\in\{-1,1\}.
\end{equation}

The inverse gnomonic map
\begin{equation}\label{eq:gnomonic}
 (X,Z)\longmapsto\frac{(X,-1,Z)}{\sqrt{1+X^2+Z^2}}
\end{equation}
takes the rectangle $D=[-u,u]\times[-v,v]$ onto $Q_{u,v}$. Its area element is
$(1+X^2+Z^2)^{-3/2}\,dX\,dZ$. Integrating this area element over the rectangle gives
\begin{equation*}\label{eq:area}
 A(u,v)=|Q_{u,v}|
       =4\arctan\frac{uv}{\sqrt{1+u^2+v^2}}.
\end{equation*}

\begin{figure}[tbp]
\centering
\begin{tikzpicture}[>=Latex,font=\small]
\def\qu{0.5898238}
\def\qv{1.913439}
\pgfmathsetmacro{\qb}{sqrt(1+\qu*\qu+\qv*\qv)}
\begin{scope}[x=2.6cm,y=2.6cm]
 \fill[blue!3] (0,0) circle (1);
 \draw[gray!60] (0,0) circle (1);
 \draw[->,gray!70] (-1.12,0)--(1.16,0) node[right] {$x$};
 \draw[->,gray!70] (0,-1.1)--(0,1.13) node[above] {$z$};
 \path[fill=blue!20,draw=blue!70!black,thick]
   plot[domain=-\qu:\qu,samples=70,variable=\t]
     ({\t/sqrt(1+\t*\t+\qv*\qv)},{\qv/sqrt(1+\t*\t+\qv*\qv)})
   -- plot[domain=\qv:-\qv,samples=100,variable=\t]
     ({\qu/sqrt(1+\qu*\qu+\t*\t)},{\t/sqrt(1+\qu*\qu+\t*\t)})
   -- plot[domain=\qu:-\qu,samples=70,variable=\t]
     ({\t/sqrt(1+\t*\t+\qv*\qv)},{-\qv/sqrt(1+\t*\t+\qv*\qv)})
   -- plot[domain=-\qv:\qv,samples=100,variable=\t]
     ({-\qu/sqrt(1+\qu*\qu+\t*\t)},{\t/sqrt(1+\qu*\qu+\t*\t)}) -- cycle;
 \draw[dashed,blue!50] (-\qu/\qb,\qv/\qb)--(\qu/\qb,-\qv/\qb);
 \draw[dashed,blue!50] (\qu/\qb,\qv/\qb)--(-\qu/\qb,-\qv/\qb);
 \foreach \sx/\sz/\lab/\anc in {1/1/{V_{++}}/south west,-1/1/{V_{-+}}/south east,1/-1/{V_{+-}}/north west,-1/-1/{V_{--}}/north east}{
  \fill[blue!70!black] (\sx*\qu/\qb,\sz*\qv/\qb) circle (0.6pt)
       node[anchor=\anc,inner sep=3pt] {$\lab$};
 }
 \fill (0,0) circle (0.7pt) node[anchor=north west,fill=blue!20,inner sep=2pt] {$c$};
 \node[align=center] at (0,-1.35) {Orthographic view from $y<0$\\on the unit sphere};
\end{scope}
\begin{scope}[shift={(7.6cm,0)},x=1.3cm,y=1.3cm]
 \fill[blue!20] (-\qu,-\qv) rectangle (\qu,\qv);
 \draw[blue!70!black,thick] (-\qu,-\qv) rectangle (\qu,\qv);
 \draw[->,gray!70] (-1.25,0)--(1.4,0) node[right] {$X$};
 \draw[->,gray!70] (0,-2.25)--(0,2.3) node[above] {$Z$};
 \fill (0,0) circle (1.3pt) node[below right] {$c$};
 \node[below left] at (-\qu,0) {$-u$};
 \node[below right] at (\qu,0) {$u$};
 \node[above left] at (0,\qv) {$v$};
 \node[below left] at (0,-\qv) {$-v$};
 \node[align=center] at (0,-2.7) {Gnomonic image\\$[-u,u]\times[-v,v]$};
\end{scope}
\draw[->,thick] (3.05,0.5)--(5.6,0.5)
 node[midway,above,align=center] {$X=x/(-y)$\\$Z=z/(-y)$};
\end{tikzpicture}
\caption{The actual quadrilateral used to exclude seven cylinders of radius
$0.9978$: $u=0.5898238$, $v=1.913439$. Its four edges are great-circle arcs. The point $c$ is
$(0,-1,0)$ and the vertices are given by \eqref{eq:vertices}.}
\label{fig:Q}
\end{figure}
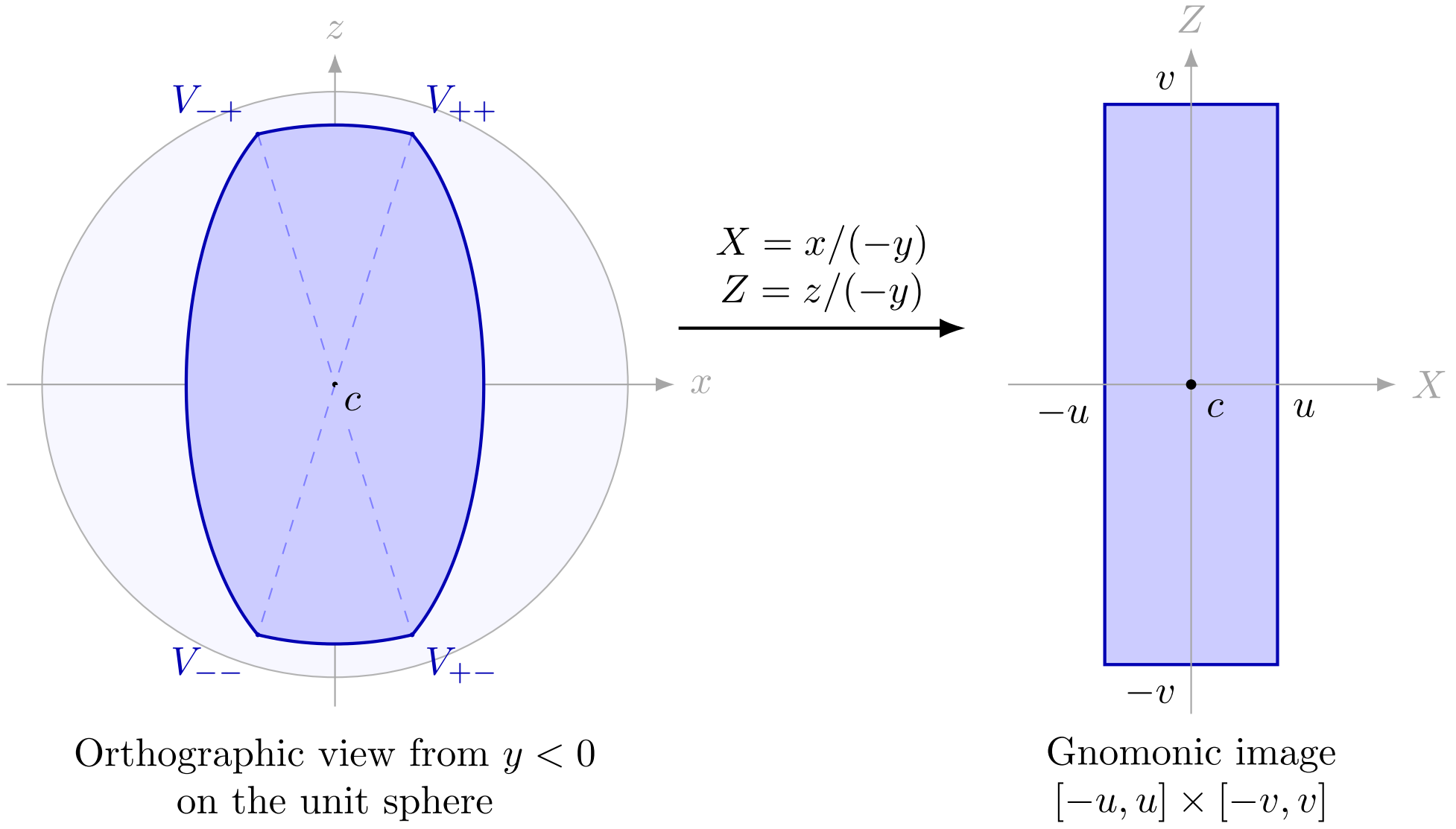

Table~\ref{tab:parameters} in the appendix gives all the rational geometric parameters. In every row all conditions (\ref{it:1})-(\ref{it:3}) are satisfied and certified, that is, there are no triple intersections, all directed differences are connected, and the uniform overlap is at most the displayed $\delta$. The rows $N=6,7,9,\ldots,19$ give the improvements stated in Theorem~\ref{thm:bounds}; the other rows record the method's performance where it does not improve the earlier bounds.

In the appendix, we also provide a brief explanation of how conditions (\ref{it:1})-(\ref{it:3}) were verified.

\section*{Data}
Verifier package and certificates for conditions (\ref{it:1})-(\ref{it:3}) for all quadrilaterals are available at \url{https://github.com/alexey-glazyrin/cylinder_verifier}.

\section*{Acknowledgments}

Alexey Glazyrin was partially supported by the NSF grants DMS-2054536, DMS-2349063 and the Simons Foundation's Travel Support for Mathematicians program.

\section*{AI statement}

ChatGPT-6 Astra was used for producing computational certificates, optimizing the parameters, and searching the literature. The construction and proofs are due to the author. The text was written and edited with the assistance of AI. The author takes full responsibility for its correctness.

\addresshere
\newpage

\appendix
\section{Parameters of quadrilaterals}
\begin{table}[h]
\centering\small
\setlength{\tabcolsep}{4.5pt}
\begin{tabular}{rrrrrr}
\toprule
$N$&Excluded $r$&$u$&$v$&$\delta$&Margin $>$\\
\midrule
4&2.686&1.12686668&5.47652733&0.11044037&$0.00031215$\\
5&1.705&0.84607699&3.09788390&0.07477855&$0.00014222$\\
6&1.252&0.68834405&2.31313051&0.04751210&$0.00049836$\\
7&0.9978&0.5898238&1.913439&0.034065&$0.000068701$\\
8&0.838&0.52318403&1.66739516&0.02867723&$0.00068596$\\
9&0.726&0.47034729&1.51585524&0.02416033&$0.0017693$\\
10&0.642&0.42802836&1.40592755&0.02040233&$0.00094367$\\
11&0.577&0.39390061&1.31982857&0.01764770&$0.00039737$\\
12&0.526&0.36618246&1.25113951&0.01609503&$0.0038743$\\
13&0.483&0.34211097&1.19208353&0.01406377&$0.00104412$\\
14&0.448&0.32199943&1.14304753&0.01295562&$0.00286214$\\
15&0.418&0.30436505&1.10015475&0.01182463&$0.00267229$\\
16&0.392&0.28876654&1.06222131&0.01072070&$0.00096453$\\
17&0.370&0.27532423&1.02949628&0.01016592&$0.00340134$\\
18&0.350&0.26290214&0.99916629&0.00923304&$0.00063446$\\
19&0.3327&0.25199299&0.97244295&0.00868393&$0.00109335$\\
20&0.318&0.24259814&0.94934581&0.00868849&$0.00780365$\\
21&0.304&0.23354194&0.92696953&0.00816695&$0.00733664$\\
22&0.291&0.22503355&0.90583596&0.00746545&$0.00350432$\\
23&0.280&0.21775576&0.88767304&0.00737071&$0.00772738$\\
24&0.269&0.21040625&0.86921865&0.00664716&$0.00189491$\\
\bottomrule
\end{tabular}
\caption{The full family of certified parameters, including the
non-improving cases omitted from Theorem~\ref{thm:bounds}. Margin means
$NA-(3N-6)\delta-4\pi$; its displayed lower bounds are rounded down.
All other decimals in this table are exact.}
\label{tab:parameters}
\end{table}

\section{No triple intersections}

We use the main lemma from \cite{BK91}. For completeness we give a short proof of this lemma here.

\begin{lem}\cite[Main Lemma]{BK91}\label{lem:BK_main}
For any three non-overlapping infinite circular cylinders of radius $r$ in $\mathbb{R}^3$, there is no point in $\mathbb{R}^3$ at distance $<\frac 2 {\sqrt{3}} r$ from the axes of all three cylinders.
\end{lem}

\begin{proof}
Assume $P$ is at distance at most $q$ from all three axes. Denote by $a_i$, $1\leq i\leq 3$, the nearest point on axis $i$ to $P$ and write $\overline{a}=1/3 (a_1+a_2+a_3)$. Since the cylinders are non-overlapping, $\norm{a_i-a_j}\geq 2r$ for $i\neq j$. Then
\[
12r^2 \leq \sum_{i< j} \norm{a_i-a_j}^2 = 3\sum \norm{a_i-\overline{a}}^2\leq 3\sum \norm{P-a_i}^2\leq 9q^2.
\]
Therefore, $q\geq \frac 2 {\sqrt{3}} r$.
\end{proof}

In order to verify that three quadrilaterals cannot have a common point, it is sufficient to show that, otherwise, there exists a point violating the conclusion of \Cref{lem:BK_main}. We use the construction from Section \ref{sect:quad} and the notation from there. Fix $\rho>0$, and for a quadrilateral $Q(C)\subset\mathbb{S}^2$ consider $\rho Q\subset\rho\mathbb{S}^2$. Then we can find the maximal distance $q$ from points in $\rho Q$ to the axis of $C$ explicitly:
\begin{equation}\label{eq:q}
q^2=\max\left\{ \rho^2 - \frac {2(1+r)\rho} {\sqrt{1+u^2}} + (1+r)^2, \frac {\rho^2 (1+u^2)}{1+u^2+v^2} - \frac {2(1+r)\rho}{\sqrt{1+u^2+v^2}} + (1+r)^2\right\}.
\end{equation}

If for three cylinders the three corresponding regions have a common point, then there is a point $P\in\rho\mathbb{S}^2$ at distance no larger than $q$ from the axes of all three cylinders. If $q < \frac 2 {\sqrt{3}} r$, this contradicts Lemma \ref{lem:BK_main}, then there is no triple intersection. The verification consists in choosing a parameter $\rho$ given $r, u, v$, and checking that $q<\frac 2 {\sqrt{3}} r$ for $q$ found from (\ref{eq:q}).

\section{Parametrizing cylinders}

We first describe how relative positions of two cylinders are parametrized. For two cylinders of radius $r$ tangent to the unit sphere with nonparallel axes $L_1$ and $L_2$ determine a common line orthogonal to both of the axes and choose a unit vector $v$ in this line in the direction from $L_1$ to $L_2$. Then $v\cdot x$ is the same for all $x\in L_1$ and its value is in $[-1-r,1+r]$. We encode it by parameter $s\in [-1,1]$ defined by $v\cdot x = -(1+r)s$. Similarly, we define parameter $t$ by $v\cdot x = (1+r) t$ for all $x\in L_2$. For the third parameter we take the angle $w$ of rotation from $L_1$ to $L_2$, $0\leq w\le\pi$. The relative position of two cylinders is completely described by a triple $(s,t,w)$ up to isometric transformations.

The distance between the lines $L_1$ and $L_2$ is $(1+r)(s+t)$ so for pairs of non-overlapping cylinders, we have $s+t\geq \frac {2r}{1+r}$. We relax this condition to $s,t\geq \frac {2r}{1+r}-1$. We reparametrize $s=\sin\theta$ and $t=\sin\phi$, and the following parameter space:
\[
(\theta, \phi, w)\in\Omega= \left[\arcsin\frac {r-1}{1+r},\frac\pi 2\right]^2\times[0,\pi].
\]
We note that this description includes all pairs of cylinders with parallel axes as well, although these pairs do not correspond to a unique triple of parameters. 

In the implementation, $\Omega$ is enclosed in the rational box $[\lambda,P/2]^2\times[0,P]$, where $P=3.141592653589794> \pi$ and the checker verifies $\lambda < \arcsin((r-1)/(r+1))$. A certificate records a finite subdivision of this rational box into smaller boxes, obtained by repeatedly bisecting a prescribed coordinate interval. For each terminal box the verifier either proves that the separation condition is impossible throughout it, or proves the required property throughout it. It also checks that the terminal boxes form a complete subdivision of the original box. Thus there is no gap between tested orientations: a box stands for a continuous family of cylinder pairs, rather than a single sample. For overlap area one may additionally exchange the cylinders; for directed differences both orders are covered.

\section{Connected differences}

For a difference $Q_1\setminus Q_2$, we assume that $Q_1$ is defined via gnomonic coordinates as in (\ref{eq:gnomonic}). Then the difference can be obtained from the four inequalities determining the sides of the second quadrilateral:
\begin{equation}\label{eq:ell}
 \ell_k(X,Z)=a_kX+b_kZ+d_k\geq0,\qquad k=1,2,3,4.
\end{equation}
Their coefficients depend on the pair parameters. Consequently $Q_1\setminus Q_2$ corresponds exactly to the union of four sets
\[
B_k=\{(X,Z)\in D:\ell_k(X,Z)<0\}.
\]
Each $B_k$ is a convex portion of a rectangle cut off by a straight line, and is therefore connected when nonempty.

The verifier checks an even stronger condition that every nonempty $B_k$ belongs to the component containing $(0,0)$. For the first step, it marks all $B_k$ that contain $(0,0)$ for the whole parameter box. Then it uses the following simple tests.

\begin{enumerate}

\item If $B_i$ is marked, $B_j$ is unmarked, the verifier checks whether $B_i$ and $B_j$ have the same common point for all parameters from the box. If they do, then it marks $B_j$.

\item For an unmarked $B_k$, list every vertex of $D$ that could belong to $B_k$ for some parameters in the box. If each such vertex is certified to belong to the marked union throughout the box, mark $B_k$. Whenever $B_k$ is nonempty, it contains a vertex of $D$, which supplies a connection to that union.

\end{enumerate}

Repeating these tests marks all four sets in every admissible terminal box of the final certificates.

\section{Pairwise intersections}

The intersection $Q_1\cap Q_2$ corresponds to the polygon
\[
D\cap\{\ell_1\geq0\}\cap\cdots\cap\{\ell_4\geq0\},
\]
where $D=[-u,u]\times[-v,v]$ and $\ell_i$ are defined by (\ref{eq:ell}). The verifier reduces the problem to rational polygons. On each quadrant of $D$, interval bounds for the coefficients give an affine function $\ell_k^+$ that bounds $\ell_k$ above throughout the parameter box. Replacing $\ell_k \geq 0$ by $\ell_k^+\geq 0$ gives a rational polygon containing every possible intersection in that quadrant.

The implementation triangulates each enclosing polygon with triangles from the same quadrant of $D$. For a triangle whose vertices are represented by rays $a,b,c$, its spherical area is
\begin{equation}\label{eq:solid}
 2\arctan
 \frac{|\det(a,b,c)|}
 {\norm a\norm b\norm c+(a\cdot b)\norm c+
       (b\cdot c)\norm a+(c\cdot a)\norm b}.
\end{equation}
All rays here lie over one quadrant of $D$, so their pairwise dot products are positive. The denominator is positive, and rational lower bounds for it give an upper bound for the non-negative ratio in (\ref{eq:solid}).
The checker then uses $\arctan x\leq x$, or the sharper bound $\arctan x\leq x-x^3/3+x^5/5$ when $0\leq x<1$.

Summing the triangle bounds over the four quadrants gives an upper bound valid for every pair in the parameter box. Every retained terminal box passes the test that this bound is at most the prescribed $\delta$.

\end{document}